\documentclass[english,11pt]{amsart}
\usepackage{amsopn,amsthm,amsfonts,amsmath,amssymb,amscd}
\usepackage{babel}
\usepackage{amssymb}
\usepackage{graphicx}

\theoremstyle{definition}
\newtheorem{Theorem}{Theorem}
\newtheorem{Lemma}[Theorem]{Lemma}

\newtheorem{Corollary}[Theorem]{Corollary}
\newtheorem{Example}[Theorem]{Example}
\newtheorem{Remark}[Theorem]{Remark}
\newtheorem*{Def}{Definition}

\newenvironment{Proof*}{{\it Proof.}}

\newcommand{\Gl}[2]{\mathop{GL}_{#1}({#2})}
\newcommand{\Mn}[2]{\mathop{M}_{#1}({#2})}

\newcommand{\nil}[1]{\iota(#1)}

\begin{document}

\title{Invertible and nilpotent matrices over antirings}

\author{David Dol\v zan, Polona Oblak}

\address{D.~Dol\v zan:~Department of Mathematics, Faculty of Mathematics
and Physics, University of Ljubljana, Jadranska 19, SI-1000 Ljubljana, Slovenia; e-mail: 
david.dolzan@fmf.uni-lj.si}
\address{P.~Oblak: Department of Mathematics, Faculty of Computer and Information Science, University of Ljubljana,
Tr\v za\v ska 25, SI-1000 Ljubljana, Slovenia; e-mail: polona.oblak@fri.uni-lj.si}

\bigskip

\begin{abstract} 
In this paper we characterize invertible matrices over an arbitrary commutative antiring $S$ with $1$ and find the structure of
$GL_{n}(S)$. We find the number of nilpotent matrices over an entire commutative finite antiring. 
We prove that every nilpotent $n \times n$ matrix over an entire antiring can be written as a sum of $\lceil \log_2 n \rceil$ square-zero matrices and also find the necessary number of square-zero summands for an arbitrary trace-zero matrix to be expressible
as their sum. 
\end{abstract}


\maketitle 

 \section{Introduction}

A \emph{semiring} is a set $S$ equipped with binary operations $+$ and $\cdot$ such that $(S,+)$ is a commutative monoid 
with identity element 0 and $(S,\cdot)$ is a monoid with identity element 1. In addition, operations $+$ and $\cdot$
are connected by distributivity and 0 annihilates $S$. A semiring is \emph{commutative} if $ab=ba$ for all $a,b \in S$.

A semiring $S$ is called an \emph{antiring} if it is zerosumfree, i.e., if the condition $a+b=0$ implies that $a=b=0$ for all $a,b \in S$.

An antiring is called \emph{entire} if $ab=0$ implies that either $a=0$ or $b=0$.

\smallskip

For example, the set of nonnegative integers with the usual operations of addition and multiplication
is a commutative entire antiring. Boolean algebras and distributive lattices are commutative (but not entire) antirings.

\smallskip

A set $\{a_1,a_2,\ldots,a_r\} \subseteq S$ of nonzero elements is called an \emph{orthogonal decomposition of 1 of length $r$} in 
$S$ if $a_1+a_2+\ldots+a_r=1$ and $a_ia_j=0$ for all $i \ne j$. 

A matrix $A \in \Mn{n}{S}$ is \emph{an orthogonal combination} of matrices $A_1, A_2,\ldots,A_r$ if there exists an orthogonal combination
$\{a_1,a_2,\ldots,a_r\}$ of 1, such that $A=\sum_{i=1}^r a_iA_i$.

Let us denote by $U(S)$ the group of all invertible elements in $S$, i.e. $U(S)=\{a \in S; ab=ba=1 \text{ for some } b \in S\}$.

\smallskip

Tan \cite{tan07} characterized the invertible matrices over a commutative antiring.
He proved that for a commutative antiring $S$, where $U(S)=\{1\}$, a matrix $A \in \Mn{n}{S}$ is invertible if and only if 
$A$ is an orthogonal combination of some $n \times n$ permutation matrices \cite[Prop. 3.7]{tan07}. 
Here, we generalize this result to an arbitrary commutative antiring (see Theorem \ref{thm:invertible}) and
we prove that $GL_{n}(S) \simeq U(S)^n \rtimes \left(S_n\right)^k$, where $k$ is the maximal length of an orthogonal decomposition of $1$.

Tan \cite{tan08} characterized the nilpotent matrices in terms of principal permanental minors and permanental adjoint matrices. 
In Section \ref{sec:nilpotent}, we give two results on nilpotent matrices. In Theorem \ref{thm:numbernilpotent}, 
we find the number of all nilpotent matrices over an entire commutative finite antiring. 
Next, we develop a result similar to \cite{wawu91} and prove that every nilpotent $n \times n$ matrix can be written as a sum of 
$\lceil \log_2 n \rceil$ square-zero matrices.  
We also find the number of square-zero matrices needed for an arbitrary trace-zero matrix to be expressible
as their sum.

\bigskip
\bigskip

\section{Invertible matrices over S}

\bigskip

In this section, we give the characterization of invertible matrices over a commutative antiring and thus 
generalize \cite[Prop. 3.7]{tan07}.

\bigskip

\begin{Theorem}\label{thm:invertible}
 If $S$ is a commutative antiring, then $A \in M_{n}(S)$ is invertible if and only if 
  $$A=D \sum\limits_{\sigma \in S_n}{a_\sigma P_\sigma}\, ,$$
 where $D$ is an invertible diagonal matrix, $P_{\sigma}$ is a permutation matrix and 
 $\sum_{\sigma \in S_n} a_\sigma=1$ is an orthogonal decomposition of $1$.

\end{Theorem}

\smallskip

\begin{proof}
Let $A=[a_{ij}]$ be an arbitrary invertible matrix. From \cite[Theorem 3.1]{tan07}, we know that
the matrices $AA^T$ and $A^TA$ are (not necessarily equal) invertible diagonal matrices.
Then $a_{ij}a_{ik}=a_{ji}a_{ki}=0$ for all $j \neq k$. Denote by $L_i$ the entry 
$(AA^T)_{ii}=\sum\limits_{k=1}^n{a_{ik}^2} \in U(S)$.
It can be easily seen that $(\sum\limits_{k=1}^n{a_{ik}})^2=L_i$,
so we know that $l_i=\sum\limits_{k=1}^n{a_{ik}}$
is invertible for each $i$.

We write $L=\prod\limits_{i=1}^n{\sum\limits_{k=1}^n{a_{ik}}}=\sum\limits_{\sigma \in S_n}{\prod\limits_{i=1}^n{a_{i\sigma(i)}}} \in U(S)$.
Let $a_\sigma=L^{-1}\prod\limits_{i=1}^n{a_{i\sigma(i)}}$ and verify that 
$\sum\limits_{\sigma \in S_n}{a_\sigma}=1$.  
Since $a_{ik}a_{jk}=0$ for $j \ne i$, we have $a_\sigma a_\tau = 0$ for $\sigma \neq \tau$ and $a_{\sigma}^2=a_{\sigma}$. 
This gives us an orthogonal decomposition of $1$ to a sum of idempotents. (Note that $a_{\sigma}$ may be $0$ for
some $\sigma \in S_n$.)
For $a_\sigma \ne 0$, the matrix $a_\sigma A$ has exactly one nonzero element in each row (and column). 
Since $a_{ij}l_i=a_{ij}^2$ for every $i$, there exists a permutation matrix $P_\sigma$ such that 
$a_\sigma A = a_\sigma D P_\sigma$, where $D={\rm Diag}(l_1,l_2,\ldots,l_n)$.
Thus, $A=(\sum\limits_{\sigma \in S_n}{a_\sigma})A=D \sum\limits_{\sigma \in S_n}{a_\sigma P_\sigma}$.

Now, let $A=D \sum\limits_{\sigma \in S_n}{a_\sigma P_\sigma}$, where $D={\rm Diag}(d_1,d_2,\ldots,d_n)$ 
is an invertible diagonal matrix, 
$P_{\sigma}$ a permutation matrix and $\sum_{\sigma \in S_n} a_\sigma=1$  an orthogonal decomposition of $1$.
Let us write 
$B=\sum\limits_{\sigma \in S_n}
 {a_\sigma  {\rm Diag}(d_{\sigma^{-1}(1)}^{-1},d_{\sigma^{-1}(2)}^{-1},\ldots,d_{\sigma^{-1}(n)}^{-1})P_\sigma^T}$.
Since $a_\sigma$ are orthogonal idempotents, one can easily verify that $AB=I$. By \cite[Lemma 2.1]{tan07}, it follows that
$BA=I$ and thus $A$ is invertible.
\end{proof}

\bigskip

\begin{Corollary}\label{thm:invertiblegroup}
 If $S$ is a commutative antiring, then the group $GL_{n}(S)$ of invertible $n \times n$ matrices over $S$
 is isomorphic to the group $U(S)^n \rtimes \left(S_n\right)^k$, where $k$ is the maximal length of an orthogonal decomposition of $1$.
\end{Corollary}

\smallskip

\begin{proof}
  Let $1=a_1+a_2+\ldots+a_k$ be an orthogonal decomposition of $1$ of the maximal length.
  If $1=b_1+b_2+\ldots+b_r$ is another orthogonal decomposition of $1$, then by multiplying these two
  equations, for each $i$ we get $\sum\limits_{i=1}^k\sum\limits_{j=1}^r a_i b_j=0$. 
  Suppose that $a_{i_1} b_j \ne 0$ for at least two $j$. Since the longest orthogonal decomposition of 1 is 
  of length $k$, it follows that for some $i_2$ all products $a_{i_2}b_j$ are equal to 0. Thus
  $a_{i_2}=a_{i_2}\sum\limits_{j=1}^r b_j=0$, which contradicts the definition of an orthogonal decomposition. So,
  denote by $\sigma(i)$ the only index such that $a_ib_{\sigma(i)} \neq 0$ and notice that
  $a_ib_{\sigma(i)}=a_i$.  
  
  Now, for any $l$, we have 
  $b_l=b_l(\sum\limits_{i \in \sigma^{-1}(l)}{a_i})=\sum\limits_{i \in \sigma^{-1}(l)}{a_i}$, so all the 
  summands of the second orthogonal sum are actually sums of some of the summands of the first sum.
  
  Since an invertible diagonal matrix is exactly a  matrix with invertible diagonal elements, we use can
  Theorem \ref{thm:invertible} to prove the corollary.
  Conjugation with a permutation matrix preserves diagonal matrices, therefore the group of invertible 
  diagonal matrices is a normal subgroup of the group of all invertible  matrices.  So, the group of 
  invertible matrices is indeed isomorphic to a semidirect product of diagonal matrices and sums of 
  permutation matrices.  (This is, of course, not a direct product unless $n = 1$.)
\end{proof}

\bigskip

\begin{Corollary}\label{thm:entireinvertible}
 If $S$ is an entire commutative antiring, then $\Gl{n}{S} \simeq U(S)^n \rtimes S_n$.
\end{Corollary}

\smallskip

\begin{proof}
 Let $S$ be an entire antiring and $1=a_1+a_2+\ldots+a_r$ an orthogonal decomposition of 1 of length $r$.
 By definition, $a_ia_j=0$ for $i \ne j$ and since $S$ is entire, it follows that either $a_i=0$ or $a_j=0$,  which implies that $r=1$.
\end{proof}

%
%
%
%
%
%

\bigskip
\bigskip

\section{Nilpotent matrices over $S$}\label{sec:nilpotent}

\bigskip

In \cite{tan08}, Tan characterized the nilpotent matrices over a commutative antiring by studying acyclic directed graphs.
Here, we use his notation and develop some further results on nilpotent matrices.

\bigskip
For a matrix $A \in \Mn{n}{S}$ we denote by $A(i,j)$ the entry in the $i$-th row and the $j$-th column of the matrix $A$.
\bigskip

\begin{Def}
For a matrix $A \in \Mn{n}{S}$ we define a directed graph (or simply a \emph{digraph}) $D(A)$ with vertices $\{1,2,\ldots,n\}$.
A pair $(i,j)$ is an edge of $D(A)$ if and only if $A(i,j)\ne 0$.
A \emph{path} in the digraph $D(A)$ (of length $k$) is a sequence of edges 
$(i_0,i_1)$, $(i_1,i_2)$, $(i_2,i_3)$,...,$(i_{k-1},i_k)$ such that 
$A(i_0,i_1) A(i_1,i_2) A(i_2,i_3)\ldots A(i_{k-1},i_k) \ne 0$.
If $i_0=i_k$, then the path is called the \emph{cycle}. An edge $(i,i)$ is called a \emph{loop}.
A digraph is called \emph{acyclic} if it does not 
contain cycles of any length.
\end{Def}

\medskip

We will assume that $S$ has no nonzero nilpotent elements.

By Tan \cite[Prop. 3.4]{tan08}, we know that $A$ is nilpotent if and only if the digraph $D(A)$ is acyclic. 
As is well known (not only in the theory of antirings), digraphs are a useful alternative way of considering nilpotent 
matrices. For example, we have the following:

\begin{Lemma}\label{thm:path}
 Suppose that $S$ is an entire commutative antiring and let $\nil{A}$ denote the index of nilpotency of a nilpotent matrix $A\in \Mn{n}{S}$. 
 Then the longest path in the digraph $D(A)$ is equal to $\nil{A}-1$. 
\end{Lemma}

\smallskip

\begin{proof}
  Since $A^{\nil{A}}=0$ and $A^{\nil{A}-1}\ne 0$, there exists a sequence of integers 
  $i_1,i_2\ldots,i_{\nil{A}-2}$ such that 
  $A^{\nil{A}-1}(i,j)=A(i,i_1) A(i_1,i_2) \ldots A(i_{\nil{A}-3},i_{\nil{A}-2}) A(i_{\nil{A}-2},j) \ne 0$.
  Therefore, the length of the longest path in $D(A)$ is greater than or equal to $\nil{A}-1$.
  
  Suppose that there exists a path in $D(A)$ of length $\nil{A}$ that contains the edges
  $(j_0,j_1)$, $(j_1,j_2)$, $(j_2,j_3)$,...,$(j_{\nil{A}-1},j_{\nil{A}})$. Since $S$ is entire, it follows that 
  $A^{\nil{A}}(j_0,j_{\nil{A}})\ne 0$, which contradicts the definition of $\nil{A}$.
\end{proof}

\bigskip

Similarly as in Corollary \ref{thm:invertiblegroup}, we would like to describe the set of nilpotent matrices over
a commutative antiring. 
Unfortunately, the set of nilpotent matrices is not closed under addition and under multiplication. 
Let us start by giving some examples of nilpotent matrices over finite antirings.

\bigskip

\begin{Example}\label{ex:idempotent}
 Let $R$ be the lattice of all idempotents of a commutative artinian ring and let $e_1, e_2,\ldots,e_m$ be the minimal idempotents in $R$.
 Thus, $R$ consists of $2^m$ elements.
 
 Since $R$ has no nonzero nilpotent elements, the nilpotent matrices over $R$ must have all the diagonal entries equal to 0.
 Thus, all the $2 \times 2$ nilpotent matrices over $R$ are of the form  
 $\left[ \begin{matrix} 
   0 & a\\
   b & 0
  \end{matrix} \right]$, 
 where $ab=0$.  
 
 Since all the elements of $R$ are the sums of minimal idempotents, it follows that 
 $a=e_{\pi(1)}+e_{\pi(2)}+\ldots+e_{\pi(k)}$ and $b=e_{\sigma(1)}+e_{\sigma(2)}+\ldots+e_{\sigma(l)}$,
 for some $0 \leq k,l \leq m$ and permutations $\pi$ and $\sigma$ of the set $\{1,2,\ldots,m\}$.
 Since $ab=0$, it follows that $\sigma(i) \notin \{\pi(1),\pi(2),\ldots,\pi(k)\}$ for $i=1,2,\ldots,l$.
 Thus, we can set $b$ to be equal to any sum of minimal idempotents that is not represented in $a$. 
 Since there are $2^{m-k}$ such sums, the number of nilpotent  $2 \times 2$ matrices over $R$ is equal to
 $\sum_{k=0}^{m} {m \choose k} 2^{m-k}=3^m$.
 \hfill$\square$
\end{Example}

\bigskip

\begin{Example}
 It is easy to find the number of nilpotent $2 \times 2$ matrices over a finite entire commutative antiring with $q$ elements.
 
 Since the diagonal entry of a nilpotent matrix must be equal to 0, and the digraph of the nilpotent matrix is acyclic,
 a nonzero nilpotent $2 \times 2$ matrix is either of the form 
 $\left[ \begin{matrix} 
  0 & a\\
  0 & 0
  \end{matrix} \right]$ or 
 $\left[ \begin{matrix} 
  0 & 0\\
  b & 0
  \end{matrix} \right]$. Since $a, b \in S -\{0\}$ are arbitrary, it follows that there are $2q-1$ nilpotent $2 \times 2$
  matrices.
  \hfill$\square$
\end{Example}
 
\bigskip

\begin{Example}
 There are exactly $6q^3-6q^2+1$ nilpotent $3 \times 3$ matrices over a finite entire commutative antiring $S$ with $q$ elements. 
 
 Namely, the nonzero nilpotent $3 \times 3$ matrices over $S$ have either $1$, $2$, or $3$ nonzero entries. 
 And thus, their digraphs are isomorphic to one of the digraphs $D_1$, $D_2$, $D_3$, $D_4$ and $D_5$ in the Figure \ref{fig:3}. 
   
        \begin{figure}[htb]
         \begin{center}
          \includegraphics[height=2.8cm, width=14.2cm]{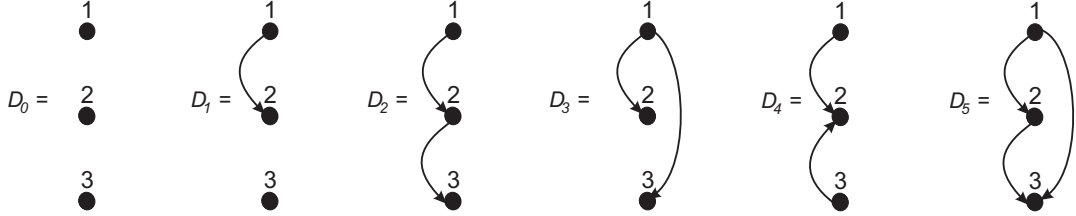}
         \end{center}
         \caption{Nonisomorphic labeled acyclic digraphs with 3 vertices}\label{fig:3}
         \end{figure}

 Note that there are $6(q-1)$ nilpotent matrices over $S$ with the digraph isomorphic to $D_1$, 
 $6(q-1)^2$ nilpotent matrices over $S$ with the digraph isomorphic to $D_2$, 
 $3(q-1)^2$ nilpotent matrices with the digraph isomorphic to $D_3$,
 $3(q-1)^2$ nilpotent matrices with the digraph isomorphic to $D_4$ and 
 $6(q-1)^3$ nilpotent matrices with the digraph isomorphic to $D_5$. Thus, there are $6q^3-6q^2$ nonzero nilpotent $3 \times 3$
 matrices over $S$.
 \hfill$\square$
\end{Example}

\bigskip

These examples give rise to a problem (and give an idea) on how to find all nilpotent matrices over a finite antiring.
The following theorem gives us the answer for a finite entire antiring.

\bigskip

\begin{Theorem}\label{thm:numbernilpotent}
 Let $S$ be a finite entire commutative antiring with $|S|=q$.  
 Then, the number of nilpotent $n\times n$ matrices over $S$ is equal to 
 $$ \sum_{\tiny{
              \begin{array}{c}
                \mu_1\geq \mu_2\geq\ldots\geq \mu_k \\
                \mu_1+\mu_2+\ldots+\mu_k=n
                \end{array}}}
         (-1)^{n-k} \frac{n!}{\mu_1!\, \mu_2!\, \ldots\, \mu_k!} q^{\frac{1}{2} \left(n^2-\sum_{i=1}^k \mu_i^2\right)}\, .
 $$
\end{Theorem}

\smallskip

\begin{proof}
 By the definition of the digraph $D(A)$, we know that every acyclic directed graph $D$ 
 corresponds to a set of nilpotent matrices $A$ (with the same pattern), such that $D=D(A)$. Namely, if $D$ 
 has $r$ edges, then nilpotent matrices $A$, such that $D(A)=D$, have $r$ nonzero entries. 

 Since $S$ is entire, it follows that the number $A_{n,r}$ of all acyclic digraphs on $n$ vertices with $r$ edges is equal to the 
 number of nilpotent $n \times n$ matrices over $S$ with exactly $r$ nonzero entries. By the main theorem of combinatorics,
 there are exactly $(q-1)^r$ such matrices.

 Thus, there are exactly $A_n(q-1)=\sum_{r=0}^{\infty}A_{n,r} (q-1)^r$ nilpotent $n\times n$ matrices.

 Let $A_n(x)$ denote the generating function for all labeled acyclic digraphs of order $n$, 
 i.e. $A_n(x)=\sum_{r=0}^{\infty}A_{n,r} x^r$. Rodionov \cite{rodionov} proved that
 %
 %
 %
  $$ A_n(x)= \sum_{m=1}^{n} (-1)^{m-1} {n \choose m} (1+x)^{m(n-m)}A_{n-m}(x)\, ,$$
 or explicitly \cite[Corollary]{rodionov},
  \begin{equation*}
  A_n(x)= \sum_{\tiny{
              \begin{array}{c}
                \mu_1\geq \mu_2\geq\ldots\geq \mu_k \\
                \mu_1+\mu_2+\ldots+\mu_k=n
                \end{array}}}
         (-1)^{n-k} \frac{n!}{\mu_1!\, \mu_2!\, \ldots\, \mu_k!} (1+x)^{\frac{1}{2} \left(n^2-\sum_{i=1}^k \mu_i^2\right)}\, .
  \end{equation*}
 Thus, the theorem holds.
\end{proof}

\bigskip

 Note that we cannot omit the condition that the antiring $S$ is entire. (See Example \ref{ex:idempotent}.)

\bigskip

\begin{Example}
 Note that the leading term of the polynomial $A_n(q-1)$ is equal to $n! \, q^{{n \choose 2}}$. Moreover, we have:
 
 \begin{tabular}{|c|l|}
 \hline
  $n$ & number of $n \times n$ nilpotent matrices over an entire antiring with $q$ elements\\
 \hline
 \hline
  $1$  &  $1$  \\
 \hline
  $2 $ & $ 2q-1 $ \\
 \hline
  $3 $ & $ 6q^3-6q^2+1 $ \\
 \hline
  $4 $ & $ 24q^6-36q^5++6q^4+8q^3+1$ \\
 \hline
  $5 $ & $ 120 q^{10}-240q^9+90q^8+60q^7-20q^6-10q^4+1$ \\
 \hline
  $6 $ & $ 720q^{15}+1800q^{14}+390q^{12}-360q^{11}-79q^9+30q^8+12q^5-1$\\
 \hline
 \end{tabular}
 
 \hfill$\square$ 
\end{Example}

\bigskip

In \cite{wawu91}, Wang and Wu characterized all matrices over a field, that can be written as a sum of two square-zero 
matrices. They showed that a matrix $T$ is a sum of two square-zero matrices if and only if it is
similar to $-T$. Fong and Sourour \cite{fosou84} showed that a matrix is a sum of two nilpotent
matrices if and only if its trace is equal to zero. Using this result, Wang and Wu showed that any 
trace-zero matrix is a sum of four square-zero matrices. 

Note that over an antiring, an arbitrary matrix with its trace equal to 0
(equivalently, all the diagonal entries of the matrix are equal to $0$) can be 
written as a sum of two nilpotent matrices 
(one being its strictly upper triangular part and the other
being its strictly lower triangular part). 
However, if a matrix over an antiring with no nonzero nilpotent elements
has a nonzero diagonal entry, then it cannot be written as a sum of 
nilpotent matrices (since its corresponding digraph contains a loop). 

\smallskip

Here we prove that every nilpotent $n \times n$ matrix can be written as a sum of 
$\lceil \log_2 n \rceil$ square-zero matrices (and that this bound is sharp). 
This implies that every trace-zero matrix over an antiring 
can be written as  a sum of at most $2 \lceil \log_2 n \rceil$ square-zero matrices. 
However, this bound is not sharp (see Example \ref{ex:sharpbound})
and we find the exact upper bound in Theorem \ref{thm:sharpbound}.

\bigskip

\begin{Def}
If $\Gamma_1$ and $\Gamma_2$ are digraphs with vertices $\{1,2,\ldots,n\}$, we 
denote by $\Gamma_1 \uplus \Gamma_2$ the digraph on vertices $\{1,2,\ldots,n\}$,
where $(i,j)$ is an edge of $\Gamma_1 \uplus \Gamma_2$ if $(i,j)$ is an
edge of $\Gamma_1$ or $\Gamma_2$.
\end{Def}
 
 \medskip
 
We can easily see that the following lemma holds. 

\begin{Lemma} \label{thm:noname}
 Let $A=A_1+A_2+\ldots+A_k$ be an $n \times n$  matrix over an antiring $S$. 
 Then $D(A)=D(A_1) \uplus D(A_2) \uplus \ldots \uplus D(A_k)$. 
\end{Lemma}

\bigskip

This enables us to prove the following theorems.

\bigskip

\begin{Theorem}\label{thm:square-zero sum}
 Let $A$ be an $n \times n$ nilpotent matrix over an entire antiring $S$. Then $A$ can be written as 
  $$A=\sum_{i=1}^{\lceil \log_2 n \rceil} B_i \, ,$$
 where $B_i \in \Mn{n}{S}$ is a square-zero matrix (i.e. $B_i^2=0$).
 
 Moreover, for every $n$ there exists a nilpotent matrix $A \in \Mn{n}{S}$ such that it 
 cannot be written as a sum of $k$ square-zero matrices, where  $k< \log_2 n$.
\end{Theorem}

\smallskip

\begin{proof}
  Without any loss of generality, we can assume that a nilpotent $A \in \Mn{n}{S}$ is a strictly uppertriangular matrix 
  (see \cite[Lemma 4.1]{tan08}). By Lemma \ref{thm:path}, all paths in the digraph corresponding to a square-zero matrix 
  are of
  length at most 1.
  
  Let $\chi(\Gamma)$ be the least number of colors needed to color the edges of a graph $\Gamma$ such that no vertex is a 
  source and a sink of two edges of the same color. Equivalently, every path in $\Gamma$ has no two incident edges of the same color.
  
  Let us denote by $\tilde{A}_n$ an arbitrary strictly uppertriangular $n \times n$ nilpotent 
  matrix with $\tilde{A}_n(i,j) \ne 0$ for $1 \leq i <j\leq n$, and let $\tilde{\Gamma}_n$ 
  be its digraph. Such a digraph $\tilde{\Gamma}_n$ is called a \emph{transitive tournament}.
  Note that $\chi(D(A)) \leq \chi(\tilde{\Gamma}_n)$.
  
  Arc colorings of some special digraphs (including transitive tournaments) were studied by
  Harner and Entringer in \cite{harent72} (and recently also by Zwonek in \cite{zwonek06}). 
  By \cite[Theorem 4]{harent72}, it follows that $\chi(\tilde{\Gamma}_n)=\lceil \log_2 n\rceil$
  and thus $\chi(D(A)) \leq \lceil \log_2 n\rceil$. 
  By Lemma \ref{thm:noname}, it follows that every nilpotent 
  $n \times n$ matrix can be written as a sum of at most 
  $\lceil \log_2 n \rceil$ square-zero matrices and $\tilde{A}_n$ cannot be written as a sum of 
  less than $\lceil \log_2 n \rceil$ square-zero matrices.
\end{proof}

\bigskip

\begin{Example}\label{ex:sharpbound}
The theorem immediately implies that every trace-zero matrix over an antiring 
can be written as a sum of at most $2 \lceil \log_2 n \rceil$ square-zero matrices. 

However, consider an arbitrary $3 \times 3$ trace-zero matrix $A$ over an antiring $S$.
 Clearly, $$A= \left[ \begin{matrix}
                0 & a & b\\
                c & 0 & d\\
                e & f & 0
             \end{matrix}
         \right]=\left[ \begin{matrix}
                0 & a & 0\\
                0 & 0 & 0\\
                0 & f & 0
             \end{matrix}
         \right]+
         \left[ \begin{matrix}
                0 & 0 & b\\
                0 & 0 & d\\
                0 & 0 & 0
             \end{matrix}
         \right]+
         \left[ \begin{matrix}
                0 & 0 & 0\\
                c & 0 & 0\\
                e & 0 & 0
             \end{matrix}
         \right]\, ,$$ 
  so $A$ can be written as a sum of $3 < 4=2 \lceil \log_2 3 \rceil$ 
  square-zero matrices.         
  \hfill$\square$
\end{Example}

\bigskip

\begin{Theorem}\label{thm:sharpbound}
 Let $A$ be an $n \times n$ trace-zero matrix over an antiring $S$ without nonzero nilpotent elements. 
 Then $A$ can be written as 
   $$A=\sum_{i=1}^{N(n)} B_i \, ,$$
 where $B_i \in \Mn{n}{S}$ is a square-zero matrix and $N(n)$ is the smallest integer, such that
 $n \leq \Big( {N(n) \atop \lceil \frac{N(n)}{2}\rceil} \Big)$.
 
 Moreover, for every $n$ there exists a trace-zero matrix $A \in \Mn{n}{S}$ such that it 
 cannot be written as a sum of $k$ square-zero matrices, where  $k< N(n)$.
\end{Theorem}

\smallskip

\begin{proof}
  Again, every path in the digraph corresponding to a square-zero matrix is of length at most 1.
  
  Recall the definition of $\chi(\Gamma)$ from the proof of the Theorem \ref{thm:square-zero sum}.  
  Let us denote by $\tilde{C}_n$ an arbitrary $n \times n$ trace-zero matrix with 
  $\tilde{C}_n(i,j) \ne 0$ for all $i \ne j$ and let $\tilde{\Delta}_n$ 
  be its digraph. Such a digraph $\tilde{\Delta}_n$ is called a \emph{complete digraph} on $n$
  vertices and it was proved in \cite[Corollary 4]{zwonek06} that 
  $\chi(\tilde{\Delta}_n)=N(n)$.
  Thus, it follows that $\chi(D(A)) \leq \chi(\tilde{\Delta}_n)=N(n)$.
  
  By Lemma \ref{thm:noname}, it follows that every trace-zero $n \times n$ matrix can be written as a sum of at most 
  $N(n)$ square-zero matrices and $\tilde{C}_n$ cannot be written as a sum of 
  less than $N(n)$ square-zero matrices.
\end{proof}

\bigskip

\begin{Remark}
 Observe that by Stirling's formula, 
 ${m \choose \lceil \frac{m}{2}\rceil} \sim \frac{2^m}{\sqrt{m}}$.  If we denote
 $\frac{2^m}{\sqrt{m}}$ by $f(m)$, then every trace-zero $n \times n$ matrix can be 
 written as a sum of approximately $f^{-1}(n)$ square-zero matrices.  Or equivalently, 
 the largest dimension of matrices, such that every trace-zero matrix can be written  as a sum of at most $n$
 square-zero matrices is approximately $\frac{2^n}{\sqrt{n}}$.  Compare this result with the
 result that the largest dimension of matrices, such that every nilpotent matrix can be written  
 as a sum of at most $n$ square-zero matrices is $2^n$.
\end{Remark}


\bigskip
\bigskip



\begin{thebibliography}{References}

  \bibitem{fosou84} 
   C.~K.~Fong, A.~R.~Sourour: \emph{Sums and products of quasi-nilpotent operators}, 
   Proc. Roy. Soc. Edinburgh A 99 (1984), 193–-200.
   
  \bibitem{harent72} 
   C.~C.~Harner, R.~C.~Entringer: \emph{Arc colorings of digraphs}, 
   J. Combinatorial Theory B 13 (1972), 219–-225.

  \bibitem{rodionov}
   V.~I.~Rodionov: \emph{On the number of labeled acyclic digraphs}, 
   Discrete Math. 105 (1992), no. 1-3, 319--321. 

  \bibitem{tan07}
   Y.~Tan: \emph{On invertible matrices over antirings}, 
   Linear Algebra Appl. 423 (2007), no. 2-3, 428--444.

  \bibitem{tan08}
   Y.~Tan: \emph{On nilpotent matrices over antirings}, 
   Linear Algebra Appl. 429 (2008), no. 5-6, 1243--1253.
  
  \bibitem{wawu91}
   J.~H.~Wang, P.~Y.~Wu: \emph{Sums of square-zero operators}, 
   Studia Math. 99 (1991), 115–-127.
   
  \bibitem{zwonek06} 
   M.~Zwonek: \emph{On arc-coloring of digraphs}, 
   Opuscula Mathematica 26 (2006), 185--195.  
\end{thebibliography}
\end{document}